\documentclass{article}
\usepackage{latexsym,amssymb,amsmath,amsfonts,pictex,enumerate,amsthm,dsfont}
\usepackage{ mathrsfs }
\usepackage{tikz}
\usetikzlibrary{decorations.markings}
\usepackage{scalerel,stackengine}
\usepackage[active]{srcltx}

\theoremstyle{definition}
\newtheorem{theorem}{Theorem}
\newtheorem{definition}{Definition}
\newtheorem{lemma}[theorem]{Lemma}

\stackMath
\newcommand\reallywidehat[1]{
\savestack{\tmpbox}{\stretchto{
  \scaleto{
    \scalerel*[\widthof{\ensuremath{#1}}]{\kern-.6pt\bigwedge\kern-.6pt}
    {\rule[-\textheight/2]{1ex}{\textheight}}
  }{\textheight}
}{0.5ex}}
\stackon[1pt]{#1}{\tmpbox}
}
\title{Bounded point evaluations in spaces of several complex variables}
\author{Stephen Deterding}
\date{Marshall University}

\begin{document}

\maketitle

\begin{abstract}
    The focus of this manuscript is existence criteria for bounded point evaluations and bounded point derivations. In the single variable setting conditions for the existence of bounded point evaluations and derivations are known for many spaces of functions; however, in the several variable case far less is known. In this manuscript, we consider Campanato spaces, which contain many other function spaces as special cases, and prove sufficient conditions for the existence of bounded point evaluations in these spaces that are analogs of known results for the single variable case. We also demonstrate that for certain cases these sufficient conditions are also necessary.
\end{abstract}

\section{Introduction}

The focus of this paper is existence criteria for bounded point evaluations and bounded point derivations in various function spaces of both one and several complex variables. Bounded point derivations were first considered in the setting of uniform rational approximation on subsets of $\mathbb C$ where they generalize the concept of non-peak points. Let $X$ be a compact subset of the complex plane, $\mathbb C$, with $x \in X$, let $R_0(X)$ denote the set of rational functions with poles off $X$ and let $R(X)$ be the closure of $R_0(X)$ in the uniform norm on $X$, which is denoted by $||\cdot||_X$.

\begin{definition}
\label{bpd-R(X)}

Let $t \in \mathbb N$. $R(X)$ admits a $t$-th order \textbf{bounded point derivation} at $x$ if there exists a constant $C$ such that

\begin{align*}
    |f^{(t)}(x)| \leq C ||f||_X
\end{align*}

\noindent for all $f \in R_0(X)$.
    
\end{definition}

\noindent A bounded point derivation of order $0$ represents evaluation and is thus called a \textbf{bounded point evaluation}. In the single variable setting bounded point evaluations and derivations have been widely studied due to connections with rational approximation and the invariant subspace problem \cite{Brennan1971, Brennan2, Conway, Wolf}. In addition, the existence of a bounded point derivation demonstrates that the functions in the space possess additional analytic structure than what would otherwise be expected. An important problem is thus to determine criteria that can be used to determine when bounded point evaluations and derivations exist. In the single variable setting, Hallstrom \cite{Hallstrom} determined an existence criterion for $R(X)$ in terms of the analytic capacity of $X$.

\begin{definition}
\label{analytic capacity}

Let $K$ be a compact subset of $\mathbb C$. The \textbf{analytic capacity} of $K$, denoted $\gamma(K)$, is

\begin{align*}
    \gamma(K) := \sup f'(\infty)
\end{align*}

where the supremum is taken over all functions $f$ that are bounded by $1$ and analytic on $\hat{\mathbb C} \setminus K$ with $f(\infty) = 0$.
    
\end{definition}

\noindent In the above definition, $f'(\infty)$ is the derivative at the point at $\infty$ on the Riemann sphere which is not necessarily $\displaystyle \lim_{z\to \infty} f'(z)$. In fact, $\displaystyle f'(\infty) = \lim_{z\to \infty} z f(z)$. 

\begin{theorem}[Hallstrom]
\label{melnikov-hallstrom condition}

Suppose $X \subseteq \mathbb C$ and $x \in X$ and let $A_n(x)$ denote the annulus $\{2^{-(n+1)} < |z-x|<2^{-n}\}$. Then $R(X)$ admits a $t$-th order bounded point derivation at $x$ if and only if the series

\begin{align}
\label{weiner-eq}
    \sum_{n=1}^\infty 2^{n(t+1)} \gamma(A_n(x) \setminus X)
\end{align}

\noindent converges.
    
\end{theorem}

\noindent A series involving capacities of the compliments of a set in concentric annuli such as \eqref{weiner-eq} is sometimes called a Wiener series due to its appearance in Wiener's criteria for a boundary point of a domain in $\mathbb R^3$ to be a regular point for the Dirichlet problem \cite[pg.130]{Wiener}. 

Bounded point evaluations and derivations have been studied for function spaces other than $R(X)$ in the single variable setting. Criteria similar to \eqref{weiner-eq} but with analytic capacity replaced with a different capacity are known to be both necessary and sufficient for the existence of bounded point evaluations and derivations in various function spaces including closures of rational functions in $L^p$ spaces \cite{Hedberg}, Lipschitz spaces \cite{Lord}, BMO spaces \cite{Deterding}, and Campanato spaces \cite{Abshire}. These examples demonstrate that in the one variable setting, existence criteria for bounded point evaluations and derivations are given in terms of a series involving a certain capacity dependent on the space of functions under consideration. Thus, in order to determine existence conditions for a given function space, one must correctly identify the associated capacity for the space. 

It is reasonable to suspect that this principle holds for the setting of several complex variables as well; however, much less is known about bounded point evaluations and derivations in spaces of several complex variables compared with the single variable case. Bounded point evaluations have been studied in $\mathbb C^d$ due to their connections to both polynomial and rational approximations \cite{ Lawrence, Range}, but little has been done to determine existence criteria similar to \eqref{weiner-eq}. The author recently determined a sufficient criterion, given in terms of a Wiener series, for the existence of bounded point evaluations for the closure of $R_0(X)$ in the $L^p$ norm, $p >2$, where $X \subseteq \mathbb C^d$ \cite{Deterding2025}; other than that, he is unaware of any other existence criteria for bounded point evaluations or derivations in the several variable setting. The purpose of this paper is to address these deficiencies by establishing existence criteria for bounded point evaluations for spaces of Campanato functions. Campanato functions generalize functions of bounded mean oscillation and include many other types of functions as special cases, such as Lipschitz spaces, thus proving existence criteria for those functions as well. Although he currently has no results for bounded point derivations, the author hopes that one day these criteria will be able to be extended to include bounded point derivations as well.

The remainder of the paper is outlined as follows. Section 2 describes Campanato functions in the setting of several complex variables, Section 3 reviews the definition and properties of lower dimensional Hausdorff content, which is the associated capacity for spaces of Campanato functions, and Section 4 establishes some preliminary results that are used in the proofs of the existence criteria. These proofs are the focus of Sections 5 and 6, with Section 5 proving the sufficiency of the existence criterion and Section 6 proving the necessity of the existence criterion for certain cases. The final section contains examples that illustrate these results.

\section{Campanato functions in several complex variables}

Campanato spaces generalize spaces of BMO functions and include both BMO spaces and Lipschitz spaces as special cases. The Campanato seminorm is defined similarly to the BMO seminorm.

\begin{definition}
\label{mean value}
    Let $X \subseteq \mathbb C^d$, $f \in L^1(X)$ and let $B$ be a ball with radius $r$. If $|B|$ denotes the volume of $B$, then \textbf{the mean value of $f$ on $B$}, denoted by $f^B$ is

\begin{align*}
    f^B = \frac{1}{|B|} \int_B |f| dV
\end{align*}

\end{definition}

\begin{definition}
\label{campanato semi-norm}

Let $1 \leq p < \infty$ and $\lambda \geq 0$ and let $X \subseteq \mathbb C^d$ be a compact set with the property that every relatively open subset of $X$ has positive volume. The \textbf{Campanato seminorm on $X$}, $||f||'_{\mathscr{L}^{p,\lambda}(X)}$, is 

\begin{align*}
    ||f||'_{\mathscr{L}^{p,\lambda}(X)} = \sup_B \left( \frac{1}{r^\lambda} \int_B |f(z)-f^B|^p dV \right)^\frac{1}{p}
\end{align*}

\noindent where $r$ is the radius of $B$ and the supremum is taken over all balls $B \subseteq X$.
    
\end{definition}

\noindent Equivalently, (See \cite[Lemma 5.6.1]{Pick} for a proof.) the Campanato seminorm can also be defined by

\begin{align*}
    ||f||'_{\mathscr{L}^{p,\lambda}(X)} = \sup_B \left( \frac{1}{r^\lambda} \inf_{c \in \mathbb C}\int_B |f(z)-c|^p dV \right)^\frac{1}{p}.
\end{align*}

\noindent where $r$ is the radius of $B$ and the supremum is taken over all balls $B \subseteq X$. In both definitions, up to a constant multiple, the supremum can be taken over cubes of side length $r$ instead of balls and we will make use of both options when convenient. When $f$ is a constant function $||f||'_{\mathscr{L}^{p,\lambda}(X)} = 0$ and thus this is not a true norm. It can be made into a norm by adding the $L^p$ norm to it.

\begin{definition}
\label{Campanato norm}

The \textbf{Campanato norm on $X$}, $||f||_{\mathscr{L}^{p,\lambda}(X)}$, is defined by

\begin{align*}
    ||f||_{\mathscr{L}^{p,\lambda}(X)} = ||f||'_{\mathscr{L}^{p,\lambda}(X)} + ||f||_{L^p(X)}
\end{align*}

\noindent where $||f||_{L^p(X)}$ denotes the $L^p$ norm on $X$.
    
\end{definition}

\begin{definition}
\label{Campanato Space}

The \textbf{Campanato space} $\mathscr{L}^{p,\lambda}(X)$ is the space of $p$-integrable functions with finite Camapanato seminorm; that is,

\begin{align*}
    \mathscr{L}^{p,\lambda}(X) = \{f \in L^p(X): ||f||'_{\mathscr{L}^{p,\lambda}(X)} < \infty\}
\end{align*}
    
\end{definition}

\noindent $\mathscr{L}^{p,\lambda}(X)$ is a Banach space endowed with the Campanato norm. Importantly, special cases of Campanato spaces coincide with other well known function spaces. For $X \subseteq \mathbb C^d$, when $2d < \lambda < 2d+p$, $\mathscr{L}^{p,\lambda}(X) = \text{ Lip}_\alpha(X)$, the space of Lipschitz continuous functions with exponent $\alpha$ on $X$ where $\alpha = \frac{\lambda-2d}{p}$ \cite[Theorem 5.7.1]{Pick}. When $\lambda < 2d$, $\mathscr{L}^{p,\lambda}(X)$ coincides with a Morrey space \cite[Theorem 5.7.1]{Pick} and with the negative Lipschitz class Lip$_\beta(X)$ where $\beta = \frac{\lambda -2d}{p}$ \cite{O'Farrell2018}. In addition, $\mathscr{L}^{p,2d}(X) = BMO(X)$ \cite[Theorem 5.8.8]{Pick} and $\mathscr{L}^{p,0}(X) = L^p(X)$ \cite[Theorem 5.5.2]{Pick}. Finally, when $\lambda > p+2d$, $\mathscr{L}^{p,\lambda}(X)$ contains only constant functions \cite[Remark 5.7.2]{Pick}. In addition, $\mathscr{L}^{p,\lambda}(X)$ coincides with other Campanato spaces, which follows from the Campanato embedding property \cite[Theorem 5.5.2]{Pick}. If $p_1,p_2, \lambda_1, \lambda_2$ are such that $1 \leq p_1, p_2 < \infty$, $0 \leq \lambda_1, \lambda_2 < \infty$, and $\frac{\lambda_1-2d}{p_1} = \frac{\lambda_2-2d}{p_2}$ then $\mathscr{L}^{p_1,\lambda_1}(X) = \mathscr{L}^{p_2,\lambda_2}(X)$.

Continuing the analogy with BMO spaces, the vanishing Campanato spaces are an analog of spaces of functions of vanishing mean oscillation (VMO) and play a similar role in the theory of Campanato spaces as spaces of VMO functions in the theory of BMO spaces and the little Lipschitz class in the theory of Lipschitz functions.

\begin{definition}
\label{vanishing campanato}

Suppose $X \subseteq \mathbb C^d$ and let $f \in \mathscr{L}^{p,\lambda}(X)$ and $\delta>0$. Then

\begin{align*}
    \Omega_f^{p,\lambda}(\delta) : = \sup_B \left\{  \left( \frac{1}{r^\lambda} \int_B |f(z)-f^B|^p dV \right)^\frac{1}{p} : \text{ radius } B \leq \delta \right\}
\end{align*}

\noindent where the supremum is taken over all balls $B \subseteq X$ with radius at most $\delta$. The \textbf{vanishing Campanato space} $V\mathscr{L}^{p,\lambda}(X)$ is the subspace of functions in $\mathscr{L}^{p,\lambda}(X)$ with the property that $\Omega_f^{p,\lambda}(\delta) \to 0$ as $\delta \to 0$.
    
\end{definition}

\noindent As in the definition of the Campanato seminorm, up to a fixed constant multiple, the supremum can be taken over cubes of side length $\delta$ instead of balls of radius $\delta$ and we will make use of both options. Like the Campanato spaces, special cases of the vanishing Campanato spaces coincide with other well-known function spaces. For $X \subseteq \mathbb C^d$, when $\lambda = 2d$, $V\mathscr{L}^{p,\lambda}(X) = VMO(X)$, the space of functions of vanishing mean oscillation on $X$ and when $2d < \lambda < 2d+p$, $V\mathscr{L}^{p,\lambda}(X) = $lip$_{\alpha}(X)$, the little Lipschitz class with exponent $\alpha$ on $X$ where $\alpha = \frac{\lambda-2d}{p}$. 

The space of analytic Camapanto functions is defined in analogy with similar spaces of analytic functions in BMO and Lipschitz spaces \cite{O'Farrell1977, Verdera}. 

\begin{definition}
\label{analytic campanato}

Let $X \subseteq \mathbb C^d$. $A_{p,\lambda}(X)$ is the space of functions in $V\mathscr{L}^{p,\lambda}(\mathbb C^d)$ that are analytic in a neighborhood of $X$.
    
\end{definition}

\noindent It follows from the inclusions mentioned above that $A_{p,\lambda}(X)$ is a space of analytic Lipschitz functions when $2d< \lambda <2d+p$ and that $A_{p,2d}(X)$ is a space of analytic BMO functions. Bounded point evaluations on these spaces are defined analogously with $R(X)$.

\begin{definition}
\label{bpe}

Let $x \in X$. $A_{p,\lambda}(X)$ admits a \textbf{bounded point evaluation} at $x$ if there exists a constant $C>0$ such that

\begin{align*}
    |f(x)| \leq C ||f||_{\mathscr{L}^{p,\lambda}(X)}
\end{align*}

\noindent for all $f \in A_{p,\lambda}(X)$.
    
\end{definition}

\section{Hausdorff Content}

In the single variable setting it is known that conditions for the existence of bounded point evaluations and derivations are given in terms of an associated capacity dependent on the function space under consideration. For example Theorem \ref{melnikov-hallstrom condition} shows that the associated capacity for $R(X)$ is analytic capacity. For a space of Campanato functions the associated capacity is a lower dimensional Hausdorff content \cite[Theorem 1]{Abshire}, which is defined as follows.

\begin{definition}
\label{Hausdorff content}
A \textbf{measure function} is a non-decreasing function $h(t)$, $t \geq 0$ such that $h(t) \to 0$ as $t \to 0$. Given a measure function $h$ and a set $E \subseteq \mathbb C^d$, the \textbf{Hausdorff content associated to $h$} is denoted $M^h(E)$ and defined as

\begin{align*}
    M^h(E) = \inf_S \sum_{Q\in S} h(r_j)
\end{align*}

\noindent where the infimum is taken over all countable coverings $S$ of $E$ by cubes $Q$ with side length $r_j$.
\end{definition}

\noindent Alternatively, up to a constant multiplicative bound, the infimum can be taken over dyadic cubes or balls of radius $r_j$.

\begin{definition}
\label{lower-dimensional hausdorff}

For $\alpha>0$ and $E \subseteq \mathbb C^d$, the \textbf{lower $\alpha$-dimensional Hausdorff content of $E$} is denoted $M_*^\alpha(E)$ and defined as

\begin{align*}
    M_*^\alpha(E) = \sup_h M^h(E)
\end{align*}

\noindent where the supremum is taken over all measure functions $h$ such that $h(t) \leq t^\alpha$ for $t \geq 0$ and $t^{-\alpha}h(t) \to 0$ as $t\to 0^+$.    
\end{definition}


Some key properties of the lower dimensional Hausdorff content are that it is a monotone set function, which means that $M_*^\alpha(E) \leq M_*^\alpha(F)$ whenever $E \subseteq F$ and that $M_*^\alpha(B) = r^\alpha$ when $B \subseteq C^d$ is a ball of radius $r$. Another important property  comes from Frostman's lemma \cite[Theorem 8.8]{Mattila} which relates lower dimensional Hausdorff content and measure.

\begin{lemma}[Frostman's lemma]
\label{Frostman}

Let $E \subseteq \mathbb C^d$ be a set with positive lower $\alpha$-dimensional Hausdorff content. Then there exists a Borel measure $\mu$ with support on $E$ such that

\begin{enumerate}
    \item $\mu(B) \leq C r^\alpha$ for all balls $B$ with radius $r$;
    \item $\displaystyle\sup_{z\in E} \frac{\mu(B(z,r))}{r^\alpha} \to 0$ as $r \to 0^+$; and
    \item $\mu(E) \geq M_*^\alpha(E)$.
\end{enumerate}
    
\end{lemma}

\section{Preliminary Lemmas}

In this section we prove various results that will be used in the proofs of the existence criteria for bounded point evaluations. We first review the Bochner-Martinelli Theorem \cite[Theorem 1.1.4]{krantz}.

\begin{definition}
The \textbf{Bochner-Martinelli kernel} is defined as

\begin{align*}
    w(\zeta,z) = \frac{(d-1)!}{(2\pi  i)^d} \frac{1}{|\zeta-z|^{2d}} \sum_{j=1}^d (\overline{\zeta_j} - \overline{z_j}) d\overline{\zeta_1} \wedge d\zeta_1 \wedge \ldots \wedge d\zeta_j \wedge \ldots \wedge d\overline{\zeta_d} \wedge d\zeta_d
\end{align*}

\noindent where the term $d\overline{\zeta_j}$ has been omitted from the wedge product. 

\end{definition}

\begin{theorem}[\textbf{Bochner-Martinelli Integral Formula}]
    Suppose that $f$ is analytic on the closure of a domain $E$ in $\mathbb C^d$ with piecewise smooth boundary $\partial E$. Then if $z\in E$,

\begin{align*}
    f(z) = \int_{\partial E} f(\zeta) w(\zeta,z)
\end{align*}

\noindent where $w(\zeta, z)$ is the Bochner-Martinelli kernel.
    
\end{theorem}

The next result appeared as Lemma 4 in \cite{Deterding2025}. Its proof is reproduced here.

\begin{lemma}
    Suppose $\zeta \neq z$. Then

\begin{align*}
    \sum_{j=1}^d \frac{\partial}{\partial\overline{\zeta_j}} \left( \frac{\overline{\zeta_j}-\overline{z_j}}{|\zeta-z|^{2d}} \right) = 0.
\end{align*}
    
\end{lemma}

\begin{proof}
To simplify calculations, we suppose that $z=0$. Since

\begin{align*}
    \frac{\partial}{\partial\overline{\zeta_j}} \left( \frac{\overline{\zeta_j}}{|\zeta|^{2d}} \right) = \frac{(\overline{\zeta_1}\zeta_1 + \ldots + \overline{\zeta_d} \zeta_d)^d-d(\overline{\zeta_1}\zeta_1 + \ldots + \overline{\zeta_d}\zeta_d)^{d-1} \overline{\zeta_j} \zeta_j}{|\zeta|^{4d}},
\end{align*}

\noindent it follows that

\begin{align*}
    \sum_{j=1}^d \frac{\partial}{\partial\overline{\zeta_j}} \left( \frac{\overline{\zeta_j}}{|\zeta|^{2d}} \right) = \frac{d(\overline{\zeta_1}\zeta_1 + \ldots + \overline{\zeta_d} \zeta_d)^d - d \sum_{j=1}^d (\overline{\zeta_1}\zeta_1 + \ldots + \overline{\zeta_d}\zeta_d)^{d-1} \overline{\zeta_j} \zeta_j}{|\zeta|^{4d}}.
\end{align*}

\noindent Using the multinomial theorem \cite[pg. 28]{Stanley} it can be shown that 

\begin{align*}
    \sum_{j=1}^d (\overline{\zeta_1}\zeta_1 + \ldots + \overline{\zeta_d}\zeta_d)^{d-1} \overline{\zeta_j} \zeta_j = (\overline{\zeta_1}\zeta_1 + \ldots + \overline{\zeta_d}\zeta_d)^d.
\end{align*}

\noindent Hence for $\zeta \neq 0$

\begin{align}
\label{vanishing derivative}
    \sum_{j=1}^d \frac{\partial}{\partial\overline{\zeta_j}} \left( \frac{\overline{\zeta_j}}{|\zeta|^{2d}} \right)=0
\end{align}

\noindent and the lemma is proved.  

\end{proof}

We now turn to the problem of estimating the Campanato norms of various functions. The next lemma shows that the integral of the mean value of $f$ is always less than or equal to the integral of $f$ itself.

\begin{lemma}
\label{mean-value-bound}
    Let $f \in L^p(X)$, $1 \leq p < \infty$ and let $q= \frac{p}{p-1}$. Let $B \subseteq \mathbb C^d$ be a ball of radius $r$ and let $f^B = \frac{1}{|B|} \int_B f dV$ denote the mean value of $f$. Then

\begin{align*}
    \left(\frac{1}{r^\lambda}\int_B |f^B|^p dV\right)^\frac{1}{p} \leq \left(\frac{1}{r^\lambda}\int_B |f|^p dV\right)^\frac{1}{p}
\end{align*}
    
\end{lemma}

\begin{proof}
    Since $|f^B|$ is a constant, it can be factored out of the integral. Thus

\begin{align*}
    \left(\frac{1}{r^\lambda}\int_B |f^B|^p dV\right)^\frac{1}{p} \leq |f^B|r^{\frac{2d-\lambda}{p}}.
\end{align*}

\noindent By H\"older's inequality, 

\begin{align*}
    |f^B| \leq \frac{1}{|B|} &\left(\int_B |f|^p dV \right)^\frac{1}{p} \left( \int_B 1 dV \right)^\frac{1}{q}\\
    &= |B|^{\frac{1}{q}-1} \left(\int_B |f|^p dV \right)^\frac{1}{p}.
\end{align*}

\noindent Hence,

\begin{align*}
    \left(\frac{1}{r^\lambda}\int_B |f^B|^p dV\right)^\frac{1}{p} &\leq r^{2d(\frac{1}{q}-1)+\frac{2d-\lambda}{p}+\frac{\lambda}{p}} \left(\frac{1}{r^{\lambda}}\int_B |f|^p dV \right)^\frac{1}{p}\\
    &=\left(\frac{1}{r^{\lambda}}\int_B |f|^p dV \right)^\frac{1}{p}
\end{align*}

\noindent which completes the proof.
    
\end{proof}

Next we consider the properties of a function $g(z)$ that will be used in the proof of Theorem \ref{main-2}. Suppose $E \subseteq \mathbb C^d$ is a set with positive lower $1+\frac{\lambda-2d}{p}$-dimensional Hausdorff content and $\mu$ is the measure supported on $E$ given by Frostman's lemma and define $g(z): = \int |z- \zeta|^{-1} d\mu(\zeta)$. Since $\mu$ has support on $E$, $g(z)$ is analytic outside $E$. In addition $g(z)$ belongs to the Campanato space $\mathscr{L}^{p,\lambda}(\mathbb C^d)$. The proof of this result is a modification of a proof of Kaufman for BMO spaces in $\mathbb C$ \cite{Kaufman} that was extended to Campanato spaces in $\mathbb C$ in \cite{Abshire}.

\begin{lemma}
\label{kaufman}
   Suppose $E \subseteq \mathbb C^d$ is a set with positive lower $\alpha$-dimensional Hausdorff content and $\mu$ is the measure supported on $E$ given by Frostman's Lemma and define $g(z) := \int |z- \zeta|^{-1} d\mu(\zeta)$. Suppose that $\lambda >0$ and let $B$ be a ball in $\mathbb C^d$ with radius $r$. Then there exists a constant $c$ and real number $C>0$ such that 

\begin{align*}
   \left( \frac{1}{r^\lambda} \int_B |g(z)-c|^p dV \right)^\frac{1}{p} \leq C
\end{align*}

\noindent and therefore $g \in \mathscr{L}^{p,\lambda}(\mathbb C^d)$
    
\end{lemma}

\begin{proof}
    Suppose $B$ has center at $w$ and let $B_* = B(w,2r)$ be the ball with the same center but twice the radius. Let 

\begin{align*}
g_1(z) = \int_{B_*} (\zeta-z)^{-1} d\mu(\zeta)
\end{align*}

\noindent and let $g_2(z) = g(z)-g_1(z)$. It follows from H\"older's inequality that

\begin{align*}
    \int_B |g_1(z)|^p dV &\leq \int_B \int_{B^*} |\zeta-z|^{-p} d\mu(\zeta) \left(\int_{B^*} d\mu(\zeta) \right)^{\frac{p}{q}} dV\\
    &\leq C r^{(1+\frac{\lambda-2d}{p})\frac{p}{q}} \int_B \int_{B^*} |\zeta-z|^{-p} d\mu(\zeta)dV,
\end{align*}

\noindent where $q = \frac{p}{p-1}$. By Fubini's Theorem,

\begin{align*}
    \int_B \int_{B^*} |\zeta-z|^{-p} d\mu(\zeta)dV &= \int_{B^*} \int_{B} |\zeta-z|^{-p} dV d\mu(\zeta)\\
    &\leq C r^{1+\frac{\lambda - 2d}{p}}r^{2d-p}. 
\end{align*}

\noindent Thus 

\begin{align*}
    \left(\int_B |g_1(z)|^p dV\right)^{\frac{1}{p}} &\leq Cr^{(1+\frac{\lambda-2d}{p})(\frac{1}{p}+\frac{1}{q})}r^{\frac{2d-p}{p}}\\
    &=C r^\frac{\lambda}{p}.
\end{align*}

\noindent Next we consider 

\begin{align*}
    \int_B |g_2(z) - g_2(w)|^p dV = \int_B \left| \int_{\mathbb C^d \setminus B^*} (\zeta-z)^{-1}-(\zeta-w)^{-1} d\mu(\zeta)\right|^p dV.
\end{align*}

\noindent Since $|\zeta-z| \leq 2|w-\zeta|$ it follows that

\begin{align*}
    \int_B |g_2(z) - g_2(w)|^p dV \leq C r^p \int_B \left( \int_{\mathbb C^d \setminus B^*} |\zeta-w|^{-2} d\mu(\zeta) \right)^p dV.
\end{align*}

\noindent To evaluate the inner integral we use a dyadic annuli estimate. Let $A_k = \{z \in \mathbb C^d: 2^kr \leq |\zeta-w| \leq 2^{k+1}r\}$. Then

\begin{align*}
    \int_{\mathbb C^d \setminus B^*} |\zeta-w|^{-2} d\mu(\zeta) \leq \sum_{k=1}^\infty (2^kr)^{-2} \mu(B(w,2^{k+1}r)) \leq Cr^{-1+\frac{\lambda-2d}{p}}.
\end{align*}

\noindent Hence,

\begin{align*}
    \int_B |g_2(z) - g_2(w)|^p dV \leq C r^p r^{2d}r^{(-1+\frac{\lambda-2d}{p})p} = C r^{\lambda}
\end{align*}

\noindent and thus

\begin{align*}
    \left(\frac{1}{r^\lambda}\int_B |g(z) - g_2(w)|^p dV\right)^{\frac{1}{p}} \leq C
\end{align*}

\noindent and the theorem is proved.
    
\end{proof}

The next result gives a bound on the size of the mean oscillation in terms of the radius of the ball.

\begin{lemma}
\label{outside-annulus-1}
    Suppose $1 \leq p < \infty$,  $0 \leq \lambda < 2d$ and $p+\lambda >2d$. Suppose $X \subseteq \mathbb C^d$ is compact, $A_n$ denotes the region $\{z\in \mathbb C^d: 2^{-(n+1)} \leq |z|\leq 2^{-n}\}$ and $\nu_n$ is the measure supported on $A_n\setminus X$ given by Frostman's lemma with $\alpha = 1+ \frac{\lambda-2d}{p}$ and let $f_n(z) = \int |z- \zeta|^{-1} d\nu_n(\zeta)$. Suppose $B$ is a ball in $\mathbb C^d$ with radius $r$. Then there exists $C>0$ such that

    \begin{align*}
        \left( \frac{1}{r^\lambda}\int_B |f_n|^p dV \right)^\frac{1}{p} \leq C r^{-(1+\frac{\lambda-2d}{p})} M_*^{1+\frac{\lambda-2d}{p}}(A_n\setminus X). 
    \end{align*}
\end{lemma}

\begin{proof}
     Since $\nu_n$ has support on $A_n \setminus X$ it follows from an application of H\"older's inequality and Fubini's theorem that

\begin{align*}
    \int_B |f_n(z)|^p dV &\leq \int_B \left( \int \frac{\nu_n(\zeta)}{|\zeta-z|^{p}}\right)\left(\int \nu_n(\zeta)\right)^\frac{p}{q}dV\\
    &\leq C  \left(M_*^{1+\frac{\lambda-2d}{p}}(A_n\setminus X)\right)^\frac{p}{q} \int \int_B \frac{1}{|\zeta-z|^{2d}} dV \nu_n,
\end{align*}

\noindent where $q = \frac{p}{p-1}$. It follows from a computation that the inner integral is bounded by a fixed constant multiple of $r^{2d-p}$ and hence

\begin{align*}
    \int_B |f_n(z)|^p dV \leq C  \left(M_*^{1+\frac{\lambda-2d}{p}}(A_n\setminus X)\right)^{\frac{p}{q}+1} r^{2d-p}.
\end{align*}

\noindent Therefore,

\begin{align*}
  \left( \frac{1}{r^\lambda} \int_B |f_n(z)|^p dV \right)^\frac{1}{p} &\leq C  M_*^{1+\frac{\lambda-2d}{p}}(A_n\setminus X)r^{-(1+\frac{\lambda-2d}{p})},
\end{align*}

which completes the proof.

\end{proof}

Finally we show that $f_n(z)$ is bounded in the $L^p$ norm provided that $1\leq p <2d$.

\begin{lemma}
\label{Cam-L^p-norm}
Let $1 \leq p <2d$ and $0 \leq \lambda \leq 2d+p$. Suppose $X \subseteq \mathbb C^d$ is compact, $A_n$ denotes the region $\{z\in \mathbb C^d: 2^{-(n+1)} \leq |z|\leq 2^{-n}\}$ and $\nu_n$ is the measure supported on $A_n\setminus X$ given by Frostman's lemma with $\alpha = 1 + \frac{\lambda-2d}{p}$ and let $f_n(z) = \int |z- \zeta|^{-1} d\nu_n(\zeta)$. Then there exists $C>0$ such that
\begin{align*}
    \left(\int_X |f_n|^p dV\right)^{\frac{1}{p}} \leq C M^{1+\frac{\lambda-2d}{p}}_*(A_n \setminus X).
\end{align*}
    
\end{lemma}

\begin{proof}
    By H\"older's inequality, 

\begin{align*}
    \int_X |f_n|^p dV 
    &\leq \int_X \left(\int \frac{\nu_n(\zeta)}{|\zeta-z|^p}\right) \left( \int\nu_n(\zeta) \right)^{\frac{p}{q}}dV\\
    &\leq C ( M^{1+\frac{\lambda-2d}{p}}_*(A_n \setminus X))^{\frac{p}{q}}\int_X \left(\int \frac{\nu_n(\zeta)}{|\zeta-z|^p}\right)dV,
\end{align*}

\noindent where $q = \frac{p}{p-1}$, and by Fubini's theorem

\begin{align*}
    \int_X \left(\int \frac{\nu_n(\zeta)}{|\zeta-z|^p}\right) dV = \int \int_X \frac{dV}{|\zeta-z|^p} \nu_n(\zeta).
\end{align*}

\noindent Let $B$ be a ball centered at $\zeta$ with radius $r$ such that $X \subseteq B$. Then

\begin{align*}
    \int_X \frac{dV}{|\zeta-z|^p} \leq \int_B \frac{dV}{|\zeta-z|^p} = C r^{2d-p}
\end{align*}

\noindent which is finite since $p<2d$. Thus

\begin{align*}
    \int_X |f_n|^p dV \leq C ( M^{1+\frac{\lambda-2d}{p}}_*(A_n \setminus X))^{\frac{p}{q}+1}
\end{align*}

\noindent and therefore

\begin{align*}
    \left(\int_X |f_n|^p dV\right)^{\frac{1}{p}} \leq C M^{1+\frac{\lambda-2d}{p}}_*(A_n \setminus X).
\end{align*}
    
\end{proof}

\section{Proof of Sufficiency}

In this section we prove a sufficient condition for $A_{p, \lambda}(X)$ to admit a bounded point evaluation when $X \subseteq \mathbb C^d$.

\begin{theorem}
\label{main}

   Suppose $1 \leq  p < \infty$,  $0 \leq \lambda \leq 2d+p$ and $p+\lambda > 2d$. Let $ X \subseteq \mathbb C^d$ be compact, $x \in X$ and $A_n(x) = \{z\in \mathbb C^d: 2^{-(n+1)} < |z-x|<2^{-n}\}$. If the sum

\begin{align*}
    \sum_{n=1}^\infty 2^{{n}} M_*^{1 + \frac{\lambda-2d}{p}}(A_n(x)\setminus X) 
\end{align*}

\noindent converges then $A_{p,\lambda}(X)$ admits a bounded point evaluation at $x$. 
    
\end{theorem}

\begin{proof}

We will show that for all $f \in A_{p,\lambda}(X)$, there exists $C$ such that $|f(x)| \leq C ||f||_{\mathscr{L}^{p,\lambda}(X)}$. Since point evaluations are unaffected by translations and dilations, we may suppose that $X$ is a subset of the unit ball in $\mathbb C^d$ and that $x = 0$. We first consider the case that $||f||_{\mathscr{L}^{p,\lambda}(X)} = 1$.

For each $n$, we can choose sets $K_n \subseteq A_n \setminus X$ so that $f$ is analytic on $A_n \setminus K_n$. Fix $n$ and let $\{Q_j\}$ be a covering of $K_n$ by dyadic cubes with no overlap except at the boundaries and let $r_j$ denote the side length of $Q_j$. Let $Q_j^*$ denote the cube with the same center as $Q_j$ but with side length of $\frac{3}{2}r_j$ and let $D_n = \cup Q_j^*$. It follows from the Bochner-Martinelli formula that 

\begin{align*}
    f(0) = \sum_{n=1}^\infty \int_{\partial D_n} f(\zeta)\omega(\zeta,0)
\end{align*}

\noindent where $\omega(\zeta,z)$ is the Bochner-Martinelli kernel. For each cube $Q_j$ construct a smooth function $\phi_j$ with support on $Q_j^*$ such that $|D^\alpha \phi_j(x)| \leq C_\alpha r_j^{-|\alpha|}$ and $\sum_j \phi_j = 1$ on a neighborhood of $\cup Q_j$. Such a construction is found in \cite[Lemma 3.1]{Harvey}. Let $\phi = 1- \sum_j \phi_j$. Then $\phi = 1$ on $\partial D_n$ and hence

\begin{align*}
    f(0) = \sum_{n=1}^\infty \int_{\partial D_n} f(\zeta)\phi(\zeta)\omega(\zeta,0).
\end{align*}

\noindent By the generalized Stokes' theorem 

\begin{align*}
    f(0) = \sum_{n=1}^\infty \int_{ D_n} \mathbf{d}[ f(\zeta)\phi(\zeta)\omega(\zeta,0)],
\end{align*}

\noindent where $\mathbf{d}$ denotes the exterior derivative. Recall that $d\overline{\zeta_1}\wedge d\zeta_1 \wedge \ldots \wedge d\overline{\zeta_d} \wedge d\zeta_d = (2i)^d dV$. After computing the exterior derivative of the Bochner-Martinelli kernel and recalling that $d\zeta_j \wedge d\zeta_k = 0$ when $j=k$, we obtain

\begin{align*}    
   f(0) &=\sum_{n=1}^\infty \int_{D_n} \sum_{k=1}^d \frac{\partial}{\partial \overline{\zeta_k}} \left( f(\zeta) \phi(\zeta) \frac{(d-1)!}{\pi^d|\zeta|^{2d}} \overline{\zeta_k} \right) dV\\
   &=\sum_{n=1}^\infty \sum_{j=1}^{M_n} \int_{Q_j^*} \sum_{k=1}^d \frac{\partial}{\partial \overline{\zeta_k}} \left( f(\zeta) \phi(\zeta) \frac{(d-1)!}{\pi^d|\zeta|^{2d}} \overline{\zeta_k} \right) dV\\
\end{align*}

\noindent where $M_n$ is the number of cubes in the covering $\{Q_j\}$. Since $\phi = 0$ on $\cup Q_j$ and $f$ is analytic on $A_n \setminus \cup Q_j$ it follows that 

\begin{align*}
    \int_{Q_j^*} \sum_{k=1}^d \frac{\partial}{\partial \overline{\zeta_k}} \left( f(\zeta) \phi(\zeta) \frac{(d-1)!}{\pi^d|\zeta|^{2d}} \overline{\zeta_k} \right)dV = \frac{ (d-1)!}{\pi^d} \int_{Q_j^*} f(\zeta) \sum_{k=1}^d \frac{\partial}{\partial \overline{\zeta_k}} \left( \phi(\zeta)\frac{\overline{\zeta_k}}{|\zeta|^{2d}} \right)dV.
\end{align*}

\noindent Next, it follows from \eqref{vanishing derivative} that 

\begin{align}
\label{rearrangement}
    \sum_{k=1}^d \frac{\partial}{\partial \overline{\zeta_k}} \left( \phi(\zeta)\frac{\overline{\zeta_k}}{|\zeta|^{2d}} \right) = \sum_{k=1}^d \frac{\partial \phi}{\partial \overline{\zeta_k}} \left(\frac{\overline{\zeta_k}}{|\zeta|^{2d}} \right)
\end{align}

\noindent and thus 

\begin{align*}
    f(0) = \frac{(d-1)!}{\pi^d} \sum_{n=1}^\infty \sum_{j=1}^{M_n} \int_{Q_j^*} f(\zeta) \sum_{k=1}^d \frac{\partial \phi_j}{\partial \overline{\zeta_k}} \left(\frac{\overline{\zeta_k}}{|\zeta|^{2d}}\right) dV.
\end{align*}

A similar calculation shows that 

\begin{align*}
    \mathbf{d}[\phi_j(\zeta)\omega(\zeta,0)] = \frac{(d-1)!}{\pi ^d}\sum_{k=1}^d \frac{\partial}{\partial \overline{\zeta_k}} \left( \phi_j(\zeta)\frac{\overline{\zeta_k}}{|\zeta|^{2d}} \right)dV.
\end{align*}

Hence it follows from \eqref{rearrangement} that

\begin{align*}
    \mathbf{d}[\phi_j(\zeta)\omega(\zeta,0)] = \frac{(d-1)!}{\pi^d} \sum_{k=1}^d \frac{\partial \phi_j}{\partial \overline{\zeta_k}} \left(\frac{\overline{\zeta_k}}{|\zeta|^{2d}}\right) dV
\end{align*}

\noindent and therefore

\begin{align*}
    \int_{Q_j^*} \frac{(d-1)!}{\pi^d} \sum_{k=1}^d \frac{\partial \phi_j}{\partial \overline{\zeta_k}} \left(\frac{\overline{\zeta_k}}{|\zeta|^{2d}}\right) dV = \int_{\partial Q_j^*} \phi_j(\zeta)\omega(\zeta,0) = 0
\end{align*}

\noindent since $\phi_j(\zeta) = 0$ on $\partial Q_j^*$. Recall that $f^{Q_j^*}$ denotes the mean value of $f$ on $Q_j^*$ Hence

\begin{align*}
    f(0) = \frac{(d-1)!}{\pi^d} \sum_{n=1}^\infty \sum_{j=1}^{M_n} \int_{Q_j^*} (f(\zeta) - f^{Q_j^*}) \sum_{k=1}^d \frac{\partial \phi_j}{\partial \overline{\zeta_k}} \left(\frac{\overline{\zeta_k}}{|\zeta|^{2d}}\right) dV.
\end{align*}

\noindent It follows from the Cauchy-Schwarz inequality that 

\begin{align*}
    |f(0)| \leq C \sum_{n=1}^\infty \sum_{j=1}^{M_n} \int_{Q_j^*} \frac{|f(\zeta)-f^{Q_j^*}|}{|\zeta|^{2d-1}} |\nabla \phi_j(\zeta)| dV.
\end{align*}

Since $\{Q_j\}$ is a covering of $K_n$, $|\zeta|> r_j$ and thus,

\begin{align*}
    |f(0)| \leq C \sum_{n=1}^\infty 2^n \sum_{j=1}^{M_n} \int_{Q_j^*} |f(\zeta)-f^{Q_j^*}| \frac{|\nabla \phi_j(\zeta)|}{r_j^{2d-2}}  dV.
\end{align*}

\noindent Hence by H\"older's inequality

\begin{align*}
    |f(0)| &\leq C \sum_{n=1}^\infty 2^{n} \sum_{j=1}^{M_n} \left(\int_{Q_j^*} |f(\zeta)-f^{Q_j^*}|^p dV\right)^{\frac{1}{p}} \left( \int_{Q_j^*} r_j^{(-2d+1)q} dV \right)^{\frac{1}{q}}\\
    &\leq C \sum_{n=1}^\infty 2^n \sum_{j=1}^{M_n} r_j^{1+\frac{\lambda-2d}{p}} \left( \frac{1}{r_j^\lambda}\int_{Q_j^*} |f(\zeta)-f^{Q_j^*}|^p dV\right)^{\frac{1}{p}}.
\end{align*}

\noindent Since $f$ belongs to $V\mathscr{L}_{p,\lambda}(\mathbb C)$, $h(t) = t^{1+\frac{\lambda-2d}{p}} \Omega_f^{p,\lambda}(\frac{3}{2}t)$ is an admissible function for $M_*^{1+\frac{\lambda-2d}{p}}(K_n)$. Thus taking the infimum over all possible covers of $K_n$ by cubes $\{Q_j\}$ gives

\begin{align*}
    |f(0)| &\leq C \sum_{n=1}^\infty 2^{n}M_*^{1+\frac{\lambda-2d}{p}}(K_n)\\
    &\leq C \sum_{n=1}^\infty 2^{n}M_*^{1+\frac{\lambda-2d}{p}}(A_n\setminus X)\\
    &\leq C.
\end{align*}

Now suppose $g \in A_{p,\lambda}(X)$ is analytic on $X$ and let $f = \frac{g}{||g||_{\mathscr{L}^{p,\lambda}(X)}}$. Then $||f||_{\mathscr{L}^{p,\lambda}(X)} = 1$ and hence $\vert f(0)\vert  \leq C$. Thus, $\vert g(0)\vert  \leq C ||g||_{\mathscr{L}^{p,\lambda}(X)}$ and there is a bounded point evaluation on $A_{p, \lambda}(X)$ at $0$.

\end{proof}

\section{Proof of Necessity}

We will show that the condition in Theorem \ref{main} is also necessary when $\lambda <2d$. When $\lambda >2d$, $A_{p,\lambda}$ is a space of analytic Lipschitz functions and thus always admits a bounded point evaluation.

\begin{theorem}
\label{main-2}

Suppose $1 \leq p < \infty$,  $0 \leq \lambda < 2d$ and $p+\lambda >2d$. Let $ X \subseteq \mathbb C^d$ be compact, $x \in X$ and $A_n(x) = \{z\in \mathbb C^d: 2^{-(n+1)} < |z-x|<2^{-n}\}$ . If $A_{p,\lambda}(X)$ admits a bounded point evaluation at $x$ then the sum

\begin{align*}
    \sum_{n=1}^\infty 2^{{n}} M_*^{1 + \frac{\lambda-2d}{p}}(A_n(x)\setminus X) 
\end{align*}

\noindent converges. 
    
\end{theorem}

\begin{proof}
    We first consider the case when $1 \leq p<2d$. By translation we may suppose that $x = 0$. Let $A_n = A_n(0)$. We will prove the contrapositive. Suppose that 

\begin{align*}
    \sum_{n=1}^\infty 2^{n} M^{1+\frac{\lambda-2d}{p}}_*(A_n \setminus X) = \infty.
\end{align*}

\noindent Then we can find a decreasing sequence $\epsilon_n \to 0$ such that 

\begin{align*}
    \sum_{n=1}^\infty \epsilon_n 2^{n} M^{1+\frac{\lambda-2d}{p}}_*(A_n \setminus X) = \infty.
\end{align*}

\noindent and $2^{n}\epsilon_nM^{1+\frac{\lambda-2d}{p}}_*(A_n\setminus X) \leq 1$ for all $n$. Frostman's lemma implies that there exists a positive measure $\nu_n$ supported on $A_n \setminus X$ such that 

\begin{enumerate}
    \item $\nu_n(B) \leq \epsilon_n r^{1+\frac{\lambda-2d}{p}}$ for all balls with radius $r$;
    \item $\sup_z \frac{\nu_n(B(z,r))}{r^{1+\frac{\lambda-2d}{p}}} \to 0$ as $r \to 0^+$; and
    \item $\int \nu_n = C\epsilon_n M_*^{1+\frac{\lambda-2d}{p}}(A_n \setminus X)$.
\end{enumerate}

Let

\begin{align*}
    f_n(z) = \int \frac{\nu_n(\zeta)}{\zeta-z}.
\end{align*}

\noindent Then $f_n(z)$ is analytic outside $A_n \setminus X$. In addition, it follows from Lemma \ref{kaufman} and Lemma \ref{Cam-L^p-norm} that $f_n \in A_{p,\lambda}(X)$. Also, $|f_n(0)| \geq  C \epsilon_n 2^n M^{1+\frac{\lambda-2d}{p}}_*(A_n \setminus X)$. For each $m \in \mathbb N$, choose $M$ so that

\begin{align*}
    1 \leq \sum_{n=m}^M \epsilon_n 2^n M^{1+\frac{\lambda-2d}{p}}_*(A_n \setminus X) \leq 2
\end{align*}

\noindent and let 

\begin{align*}
    g_m(z) = \sum_{n=m}^M f_n(z).
\end{align*}

\noindent Thus $|g_m(0)|$ is bounded below by a non-zero constant for all $m$. We wish to show that $||g_m||_{\mathscr{L}^{p,\lambda}(X)} \to 0$ as $m \to \infty$. Let $B$ be a ball of radius $r$ and choose $m \in \mathbb N$ such that $r> 2^{-m}$. Recall that $f_n^B$ denotes the mean value of $f_n$ on $B$. It follows from Lemma \ref{mean-value-bound} and Lemma \ref{outside-annulus-1} that for $n \geq m$

\begin{align}
\label{frostman-eq}
    \left( \frac{1}{r^{\lambda}}\int_B |f_n - f_n^B|^p dV \right)^\frac{1}{p} \leq C 2^{n(1+\frac{\lambda-2d}{p})} \epsilon_n M_*^{1+\frac{\lambda-2d}{p}}(A_n\setminus X). 
\end{align}

\noindent Note that the difference in the particular version of Frostman's lemma used in this proof and the general version stated in Lemma \ref{Frostman} explains the presence of the term $\epsilon_n$ in the right hand side of \eqref{frostman-eq}. It follows from the definition of $g_m$ and the Minkowski inequality that

\begin{align*}
    \left( \frac{1}{r^{\lambda}}\int_B |g_m - g_m^B|^p dV \right)^\frac{1}{p} \leq \sum_{n=m}^M\left( \frac{1}{r^{\lambda}}\int_B |f_n - f_n^B|^p dV \right)^\frac{1}{p}
\end{align*}

\noindent and hence

\begin{align*}
    \left( \frac{1}{r^{\lambda}}\int_B |g_m - g_m^B|^p dV \right)^\frac{1}{p} &\leq C \sum_{n=m}^M 2^{n(1+\frac{\lambda-2d}{p})} \epsilon_n M_*^{1+\frac{\lambda-2d}{p}}(A_n\setminus X)\\
    &\leq C 2^{m(\frac{\lambda-2d}{p})}\sum_{n=m}^M\epsilon_n 2^{n} M^{1+\frac{\lambda-2d}{p}}_*(A_n \setminus X). 
\end{align*}

\noindent Since $\lambda < 2d $ this last line tends to $0$ as $m \to \infty$. In addition, since $1 \leq p < 2d$ it follows from Lemma \ref{Cam-L^p-norm} that 

\begin{align*}
 ||g_m||_{L^p(X)} &\leq C \sum_{n=m}^M\epsilon_n M^{1+\frac{\lambda-2d}{p}}_*(A_n \setminus X)\\
 &\leq C 2^{-m}\sum_{n=m}^M\epsilon_n 2^{n} M^{1+\frac{\lambda-2d}{p}}_*(A_n \setminus X),
\end{align*}

\noindent which also tends to $0$ as $m \to \infty$. Hence $||g_m||_{\mathscr{L}^{p,\lambda}(X)} \to 0$ as $m \to \infty$ but $|g_m(0)|$ is bounded below by a non-zero constant for all $m$, which shows that $A_{p,\lambda}(X)$ does not admit a bounded point evaluation at $0$. This proves Theorem \ref{main-2} for the case that $1 \leq p < 2d$.

For the case that $p \geq 2d$ we use the coincidence of Campanato spaces. Recall that if $p_1,p_2, \lambda_1, \lambda_2$ are such that $1 \leq p_1, p_2 < \infty$, $0 \leq \lambda_1, \lambda_2 < \infty$, and $\frac{\lambda_1-2d}{p_1} = \frac{\lambda_2-2d}{p_2}$ then $\mathscr{L}^{p_1,\lambda_1}(X) = \mathscr{L}^{p_2,\lambda_2}(X)$. If $p \geq 2d$ and $\lambda < 2d$, choose $\lambda_1$ such that $2d+\frac{2d(\lambda-2d)}{p} < \lambda_1 \leq 2d+\frac{\lambda-2d}{p}$ and let $p_1 = p \left(\frac{2d-\lambda_1}{2d-\lambda} \right)$. Then $1 \leq p_1 < 2d$ and $\mathscr{L}^{p_1,\lambda_1}(X) = \mathscr{L}^{p,\lambda}(X)$, so the proof of Theorem \ref{main-2} for the pair $p$ and $\lambda$ follows from the proof with the pair $p_1$ and $\lambda_1$.

\end{proof}

\section{Examples}

In this section we present two examples that illustrate Theorems \ref{main} and \ref{main-2}. Let $1 \leq p < \infty$ and $0 \leq \lambda \leq 2d+p$ be such that $p+\lambda > 2d$. The first example is of a set $X \subseteq \mathbb C^d$ such that $A_{p,\lambda}(X)$ admits a bounded point evaluation at $0$, which is built by generalizing the construction of a roadrunner set. Let $A_n = \{z: 2^{-(n+1)} < |z| < 2^{-n}\}$ and remove from each $A_n$ a ball $B_n$ with radius $r_n = 4^{-n(\frac{p}{p+\lambda-2d})}$ chosen so that it is entirely contained in $A_n$. Let $B$ denote the closed unit ball centered at $0$ and let $X = B \setminus \bigcup_{n=1}^\infty B_n$. Since the lower $\alpha$-dimensional Hausdorff content of a ball with radius $r$ is $r^\alpha$, it follows that

\begin{align*}
    \sum_{n=1}^\infty 2^n M_*^{1+\frac{\lambda-2d}{p}}(A_n \setminus X) = \sum_{n=1}^\infty 2^{-n} = 1
\end{align*}

\noindent and thus it follows from Theorem \ref{main} that $A_{p,\lambda}(X)$ admits a bounded point evaluation at $0$.

Suppose $\lambda_0 < 2d$. The second example uses a similar construction to create a set $X \subseteq \mathbb C^d$ such that $A_{p,\lambda}(X)$ admits a bounded point evaluation at $0$ when $\lambda > \lambda_0$ but does not admit a bounded point evaluation at $0$ when $\lambda \leq \lambda_0$. To construct this set, we again start with $B$, the closed unit ball centered at $0$ but this time we remove from each $A_n$ a ball $B_n$ with radius $r_n = 2^{-n(\frac{p}{p+\lambda_0-2d})}$ and let $X = B \setminus \bigcup_{n=1}^\infty B_n$. Then

\begin{align*}
    \sum_{n=1}^\infty 2^n M_*^{1+\frac{\lambda-2d}{p}}(A_n \setminus X) = \sum_{n=1}^\infty 2^{n(\frac{\lambda_0-\lambda}{p+\lambda_0 -2d})}.
\end{align*}

When $\lambda > \lambda_0$ this is a convergent geometric series and therefore $A_{p,\lambda}(X)$ admits a bounded point evaluation at $0$ by Theorem \ref{main}; however, when $\lambda \leq \lambda_0$ this series diverges and thus $A_{p,\lambda}(X)$ does not admit a bounded point evaluation at $0$ by Theorem \ref{main-2}.

\bibliographystyle{abbrv}
\bibliography{biblio}

\end{document}